%Beginning Formats                                        

\magnification=1200
%\magnification=1000
%\magnification=\magstep1
%\vsize=20truecm
%\voffset=0.50truein
\settabs 18 \columns
%\hoffset=3.75truecm
 
%\hsize=14truecm

%\nopagenumbers
\baselineskip=17 pt
%\ifnum\pageno=1
%\topinsert \vskip 1.00 in
%\endinsert
%\vsize=7.5in
%\fi

\def\overleftarrow#1{\vbox{\ialign{##\crcr
  \leftarrowfill\crcr\noalign{\kern-1pt\nointerlineskip}
  $\hfil\displaystyle{#1}\hfil$\crcr}}}
\def\overrightarrow#1{\vbox{\ialign{##\crcr
  \leftarrowfill\crcr\noalign{\kern-1pt\nointerlineskip}
  $\hfil\displaystyle{#1}\hfil$\crcr}}}

%Ecriture des corps de nombres

\def\Cit{\hbox{\it l\hskip -5.5pt C\/}}

\def\Crm{\hskip0.5mm \hbox{\rm l\hskip -5.5pt C\/}}

\def\s{\smallskip}

\def\b{\bigskip}
\def\bb{\bigskip\bigskip}
\def\bbb{\bigskip\bigskip\bigskip}

\def\sqr#1#2{{\vcenter{\vbox{\hrule height.#2pt
 \hbox{\vrule width.#2pt height#1pt \kern#1pt
 \vrule width.#2pt} \hrule height.#2pt}}}}

\def\operp{\hbox{${\kern+.25em{\bigcirc}
\kern-.85em\bot\kern+.85em\kern-.25em}$}}
%\def\gapprox
%{\hbox{$
%\smash{lower0.5ex\hbox{$\scriptstyle>$}} \atop
%\smash{raise0.3ex\hbox{$\scriptstyle \sim$}}
%$}}
%\def\lapprox
%{\hbox{$
%\smash{lower0.5ex\hbox{$\scriptstyle<$}} \atop
%\smash{raise0.3ex\hbox{$\scriptstyle \sim$}}
%$}}
\def\lsim{\;\raise0.3ex\hbox{$<$\kern-0.75em\raise-1.1ex\hbox{$\sim$}}\;}
\def\gsim{\;\raise0.3ex\hbox{$>$\kern-0.75em\raise-1.1ex\hbox{$\sim$}}\;}
\def\no{\noindent}
\def\r{\rightline}
\def\ce{\centerline}
\def\ve{\vfill\eject}
\def\rdots{\mathinner{\mkern1mu\raise1pt\vbox{\kern7pt\hbox{.}}\mkern2mu
 \raise4pt\hbox{.}\mkern2mu\raise7pt\hbox{.}\mkern1mu}}

\def\e e{$e^+ e^-$ }

%End of Beginning Formats
%Beginning of Letter Heading
%\rightline {UCLA/95/TEP/32}

\def\today{\ifcase\month\or January\or February\or March\or April\or
May\or June\or July\or August\or September\or October\or November\or
December\fi
\space\number\day, \number\year}

\r \today 

\bbb
\ce {\bf EXACT DEFORMATIONS OF QUANTUM GROUPS;}\ce {\bf   
 APPLICATIONS TO THE AFFINE CASE}\b
\ce {C. Fronsdal}
\s
\ce {\it Department of Physics}
\ce {\it University of California, Los Angeles, CA 90095-1547}
\bb

\no {\bf Abstract.} This paper continues our investigation of a class of
generalized quantum groups. The ``standard" R-matrix was shown to be the
unique solution of a very simple, linear recursion relation and the
classical limit was obtained in the case of quantized Kac-Moody algebras
of finite type. Here the standard R-matrix for  generalized quantum groups
is first examined in the case of quantized affine Kac-Moody algebras. 
The  classical limit yields the standard affine r-matrices of Belavin and 
Drinfeld. Then, turning to the general case, we study the exact
deformations of the standard R-matrix and the  associated Hopf algebras.
They are described as a generalized twist, $ R_\epsilon = (F^t)^{-1}RF$,
where $R$ is the standard R-matrix and $F$ (a power series in  the
deformation parameter $\epsilon$) is the solution of a linear recursion
relation of the same type as that which determines
$R$.   
   Specializing again, to the case of quantized, affine Kac-Moody
algebras, and taking the  classical limit of these esoteric quantum
groups, one re-discovers the esoteric affine r-matrices of Belavin and
Drinfeld, including the elliptic ones.  The formulas obtained here are
easier to use than the original ones, and the structure of the space of
classical r-matrices (for simple Lie algebras)  is more tranparent. In
addition,  the r-matrices obtained here are more general in that they are
defined on the central extension of the loop groups.

\ve

\no {\bf 1. Introduction.}

It is now known that most Lie bi-algebras can be quantized, but there is 
not yet a workable universal construction of the corresponding quantum
groups. For simple Lie algebras, and for the affine Kac-Moody algebras,
Belavin and Drinfeld [BD] gave a  complete classification of the 
r-matrices; and thus, of the associated coboundary bi-algebras. They are 
  the most interesting Lie bialgebras and their quantized versions is one
of  the subjects of this paper.  
 
Our strategy is to construct the universal R-matrices for a class of 
generalized standard quantum groups (those that ``commute with Cartan"),
and to  discover the others by means of deformation theory. In a previous
report we have calculated all deformations of a certain type, up to first
order in the  deformation parameter. Here 
  exact deformations are obtained, to all orders in the deformation
parameter. Unfortunately, we still do not have a good characterization of
a category within which the deformations should be sought; the type of
deformations examined is thus somewhat {\it ad hoc} (see Eq.(1.14) below);
nevertheless, all the  trigonometric r-matrices are recovered in the
classical limit, with their central extensions. In addition, the elliptic
r-matrices turn up  as a special case; these enigmatic objects thus find
their natural place.

The standard R-matrix is constructed by a method that has already been
put  to some effect in a previous paper. (The same method was used by
Lusztig [L]  in a  more special context.) The present  work encompasses
all the simple quantum groups (including the multiparameter or twisted
versions constructed by Reshetikhin [R]) and   the quantized Kac-Moody
algebras, as well as other coboundary Hopf algebras that have nothing to
do with Lie groups and that have not yet been investigated in detail. In
this paper we specialize to the case of affine Kac-Moody algebras, except
that Theorem 5.1 and Proposition 5.2 hold in the general case. 

The affine Kac-Moody algebras are central extensions of loop algebras. The
r-matrices obtained here    are defined on the extended loop algebras and
are therefore in this  respect more general than those of Belavin and
Drinfeld.
 
The standard R-matrix is a formal series
$$ R = {\rm exp}\bigl(\varphi^{ab}H_a \otimes H_b\bigr)\Bigl(1 +
\sum_{n=1}^\infty t_n\Bigr),
\quad  t_n = \sum t_{(\alpha)}^{(\alpha')}\,e_{-\alpha_1}\ldots
e_{-\alpha_n}
\otimes e_{\alpha_1'}\ldots e_{\alpha_n'}.\eqno(1.1)
$$
 Here
$t_{(\alpha)}^{(\alpha')}$ are complex coefficients and
$H_a, e_{\pm\alpha}$ are Chevalley-Drinfeld generators. 

More precisely: 

\vskip.50cm

 \no {\bf Definition.}  Let $M,N$ be two countable sets,
$\varphi,H$  two maps,
$$
\eqalign{& \varphi :~~M\otimes M \rightarrow \Crm~, \cr & H :~~M\otimes N
\rightarrow
\Crm~, \cr} \quad
\eqalign{a,b &\rightarrow \varphi^{ab}~, \cr a,\beta & \rightarrow
H_a(\beta)~. \cr}
\eqno(1.2)
$$
\vskip.5cm
\no Let ${\cal{A}}$ or ${\cal{A}}(\varphi,H)$ be the universal,
associative, unital  algebra over \Crm \enskip with generators
$\{H_a\}\, a\in M,~ \{e_{\pm \alpha}\}\,\alpha \in N$, and relations
$$
\eqalignno{&[H_a,H_b]=0~, \quad  [H_a,e_{\pm\beta}] = \pm
H_a(\beta)e_{\pm\beta}~, & (1.3)
\cr &[e_\alpha,e_{-\beta}]=\delta^\beta_\alpha
\bigl(e^{\varphi(\alpha,\cdot)}-e^{-\varphi(\cdot,\alpha)}\bigr)~, & 
(1.4) \cr}
$$
\no with $\varphi(\alpha,\cdot)=\varphi^{ab}H_a(\alpha)H_b,~
\varphi(\cdot,\alpha)=\varphi^{ab}H_aH_b(\alpha)$ and $ e^{\varphi(\alpha,
\cdot) +
\varphi(\cdot,\alpha)} \neq 1,  \alpha \in N
$.    The free subalgebra generated by
$\{e_\alpha\}~\alpha\in N$ (resp.
$\{e_{-\alpha}\}~ \alpha\in N$) will be  denoted
${\cal{A}}^+$ (resp.
${\cal{A}}^-$). The subalgebra ${\cal{A}}^0$ generated by $\{H_a\}\, a\in
M$ will be  called the Cartan subalgebra.
\vskip.5cm

The sum in (1.1) is over all $\alpha_i \in N$ and all permutations
$(\alpha')$ of the set
$(\alpha)$. We set 
$$  t_1 = \sum_{\alpha \in N} e_{-\alpha} \otimes e_\alpha~.\eqno(1.5)
$$  If the parameters of ${\cal A}$ are in general position, then all the
other coefficients are determined by the requirement that $R$ satisfy the
Yang-Baxter relation
$$  R_{12}R_{13}R_{23} = R_{23}R_{13}R_{12}.
\eqno(1.6)
$$  For special values of the parameters the R-matrix is still uniquely 
determined by the Yang-Baxter relation, but now it has to be defined on a
quotient ${\cal A}/I$ where   $I \in {\cal A}$ is a suitable ideal [F].
Quantized Kac-Moody algebras are characterized by the property that, for
each pair $(\alpha,\beta)$, there is a positive integer 
$k = k_{\alpha\beta}$ such that the following relation holds
$$  e^{\varphi(\alpha,\beta) + \varphi(\beta,\alpha) +
(k-1)\varphi(\alpha,\alpha)}
 =1.\eqno(1.7)
$$  In this case the ideal $I$ is generated by the Serre relations
$$  0 = \sum^k_{m=0} Q^k_m~(e_\alpha)^me_\beta\, (e_\alpha)^{k-m}~,
\eqno(1.8)
$$  with coefficients
 $$
 Q^k_m=(-)^m e^{m\varphi(\alpha,\beta)}~q^{m(m-1)/2}
\left(\matrix{k\cr m\cr}\right)_q~,\quad q := e^{\varphi(\alpha,\alpha)}.
\eqno(1.9)
$$ 
\ve We suppose Card $N$ and Card $M$ finite and interpret $A = 1-k$ as the
generalized Cartan matrix of a Kac-Moody algebra; this allows us to extend
the terminology of the classification of Kac-Moody algebras to quantized
Kac-Moody algebras. 

Let ${\cal A}'_{cl}$ be the algebra obtained from
${\cal A}'$ when the relations (1.4) are replaced by
$$  [e_\alpha,e_{-\beta}] = \delta_\alpha^\beta\,\bigl(
\varphi(\alpha,\cdot) +
\varphi(\cdot,\alpha)\bigr).\eqno(1.10)
$$
 If ${\cal A}_{cl}'$ is a Kac-Moody algebra of finite type, resp. affine
type, then we  may say that ${\cal A}'$ is a quantized Kac-Moody algebra
of finite type, resp. affine type. But because ${\cal A}'$ cannot be
recovered from ${\cal A}_{cl}'$ an autonomous definition is preferable.
\vskip.50cm
\no {\bf Definition 1.2.} Let ${\cal A'}$ be as above; that is, the
quotient of an algebra
${\cal A}$ as per Definition 1.1, with parameters satisfying (1.7), by the
ideal  generated by the Serre relations (1.8). We shall say that ${\cal
A}'$  
 is a quantized Kac-Moody algebra of finite type if (i) Card $M$ = \break
Card $N = l <\infty$,  and (ii) the (symmetrizable) generalized Cartan
matrix
$$ A = 1 - k, \quad k_{\alpha\beta} = \bigl(\varphi(\alpha,\beta) +
\varphi(\beta,\alpha)\bigr) /\varphi(\alpha,\alpha)\eqno(1.11)
$$  is positive definite with $A_{\alpha\beta} \in \{0,-1,\ldots\}, \,\,
\alpha \neq \beta$.
\footnote * {We are here assuming that the relation (1.7) is solved by the
exponent taking the value zero.} We shall say that
${\cal A}'$  
 is a quantized Kac-Moody algebra of affine type if (i) Card $M$ = 1 + 
Card $N   <
\infty$, and (ii) the  generalized Cartan matrix
 is positive semi-definite with $A_{\alpha\beta} \in \{0,-1,\ldots\}, \,\,
\alpha \neq
\beta$ and
  all its principal minors are  positive definite.

\vskip.50cm

A quantized affine Kac-Moody algebra can be described as follows.  Let
$\hat {\cal A}'$ be as above, with parameters satisfying (1.7) and Serre
relations  (1.8), with root generators $\{e_{\pm\alpha}\}\,\alpha =
0,\ldots,l $ and Cartan generators 
$H_1,\ldots ,H_l,\,c,\,d$, such that the subset that consists of $\{e_{\pm
\alpha}\}\,
\alpha \neq 0$ and $H_1, \ldots , H_l$  
 generates a subalgebra ${\cal A}'$ that is a quantized Kac-Moody algebra
of finite type. Let $\hat \varphi$ refer to $\hat {\cal A}'$ and $\varphi$
to ${ \cal A}'$,  and suppose that
$$
\hat \varphi = \varphi + u \, c\otimes d + (1-u)\, d \otimes c,\quad  [d,
e_{\pm \alpha}] = \pm \delta_\alpha^0 \,e_{\pm 0}.
\eqno(1.12)
$$ with some $ u \in \Cit$. 
 Suppose that $c$ is central and that the extra root defined by
$ [H_a,e_0] = H_a(0)e_\alpha
$ is such as to make the generalized Cartan matrix of $\hat {\cal A}'$ 
positive semi-definite  with all its principal minors positive. Then $\hat
{\cal A}'$ is a  quantized affine Kac-Moody algebra.
 
The classical r-matrix associated with $R$ is defined after a rescaling of
the generators~-~ Eq.(2.4)~-~ as the  coefficient of $\hbar$ in  the
expansion $R = 1 + \hbar r + o(\hbar^2)$; it satisfies the classical
Yang-Baxter relation.  Note that $r + r^t \neq 0$; the antisymmetric part
of $r$ satisfies the modified classical Yang-Baxter relation.  We 
calculate this classical r-matrix, dealing separately  with the following
cases: First the unextended  loop algebras, untwisted in Section 2,
twisted in Section 3; then, in Section 4, the  full Kac-Moody algebras.
 
Our second subject is the calculation of exact deformations of the standard
R-matrix, satisfying  the Yang-Baxter relation, in the wider context of
the bialgebras 
${\cal A}$ and ${\cal A}' = {\cal A}/I$ described above. We set 
$$ R_\epsilon = R + \epsilon R_1 + o(\epsilon^2),
$$ and suppose that  $R_1$ is driven by a term of the type
$$ Se_{-\rho} \otimes e_\sigma + S'e_\sigma \otimes e_{-\rho},\,\,\,\,\,
S,S' \in {\cal A}'^0.
\eqno(1.13)
$$ Such deformations exist under certain conditions on the parameters;
then $S$ and $S'$ and the remaining terms in
$R_\epsilon$  
  (a formal power series in $\epsilon$ with constant term $R$) are
determined by the Yang-Baxter relation. An exact formula (to all orders in
$\epsilon$) in closed form  for $R_\epsilon$  is obtained for the case of
elementary deformations, when
$R_1$ is a single term  of the type (1.13). In the general case of 
compound deformations, when (1.13) is replaced by a sum of terms of the
same type,  we obtain exact deformations  in the form of a generalized
twist.

Let $R$ be the R-matrix of a coboundary Hopf algebra ${\cal A}'$, and $F
\in {\cal A}' \otimes {\cal A}'$, invertible. Then
$$
\tilde R := (F^t)^{-1}RF\eqno(1.14)
$$ satisfies the Yang-Baxter relation if $F$ satisfies the following
relation
 
$$
\bigl((1 \otimes \Delta_{21}) F\bigr)F_{12} = \bigl( (\Delta_{13} \otimes
1) F\bigr) F_{31}.\eqno(1.15)
$$ (See Theorem 5.1 for the complete statement.) Though it is not quite
germaine to our discussion,  it may be worth while to point out that, if
$R$ is unitary, then so is $\tilde R$; the formula (1.14) therefore yields
a large family of (mostly) new unitary R-matrices. 

Applying this to our context, we find that the  relation (1.15)is
equivalent to a simple, linear recursion relation that can be reduced to
the same form as the  recursion relation that determines the coefficients
in the expansion of $R$.
  It has a unique solution that can be expressed directly in terms of the
coefficients in (1.1).  Just as in the standard case, this leads to a
simple equation for the  classical r-matrix, from which the  latter is
determined to all orders.  

In Section 6 we specialize to the case of quantized, affine Kac-Moody
algebras   and  take the classical limit, to recover the esoteric, affine
r-matrices of the simple Lie algebras, with  their central extensions. The
result agrees with that of Belavin and Drinfeld, except that they  did not
include the central extension. The formulas obtained in this paper are
more transparent and simpler to use.

Finally, in Section 7, we deal with a very special case, to discover that
the elliptic r-matrices of 
$sl(N)$ form a special case among the deformed, trigonometric r-matrices.
The universal R-matrix is  expressed as an infinite product. It is shown,
in the particular case of the elliptic R-matrix for 
$sl(2)$ in the fundamental representation, that this infinite product is
both convergent and  of practical utility; it reduces to the
representation of elliptic functions  in terms of infinite products, and
the result is in perfect agreement with Baxter [B]. 
\ve

\vskip.5cm

\ve

\no {\bf 2. Untwisted loop algebras.}

Consider a quantized affine Kac-Moody algebra $\hat{\cal A}'$, with
generators $e_{\pm 0},
\ldots ,e_{\pm l}$ and $H_1,\ldots ,H_l, c, d$.   
\vskip.5cm
\no {\bf Definition 2.} Positive root vectors are elements in ${\cal
A}'^+$ defined  recursively.    (a) The generators $e_\alpha$ are positive
root vectors. (b) If 
$E_i,E_j$ are positive root vectors then so is
$$ (1-x)^{-1}(E_iE_j - e^{\varphi(i,j)}E_jE_i), \quad x = e^{-\varphi(i,j)
- \varphi(j,i)} \neq 1,
$$ (c) All positive root vectors $E_i$ are obtained in this way from the
generators. Negative root vectors are in ${\cal A }'^-$ and are defined
analogously.
\vskip.5cm It is easy to verify that
$$ [E_i,E_{-i}] = e^{\varphi(i,\cdot )} - e^{\varphi(\cdot,i)}.
$$ 

Let $\{E_i\}\,i = 1,\ldots ,n,+$ be the  positive root vectors, labelled
in such a way that
 $$ [e_\alpha,E_+] = 0 = [e_{-\alpha},E_-],\eqno(2.1)
$$  and
\vskip-4mm 
$$  [E_i, E_-]    \in {\cal A }'^0 \cdot {\cal A }'^-, 
\quad [e_{-\alpha}, E_+]   \in {\cal A}'^0 \cdot {\cal A }'^+.\eqno(2.2)
$$  Then we may refer to $E_+$ as a highest root vector.

Suppose that the extra root $H_a(0) = H_a(E_-)$, and  
 pass to the associated untwisted loop algebra ~$ \Cit 
[\lambda,\lambda^{-1}] \otimes {\cal A}'$ by substituting
$$
\hat\varphi \rightarrow \varphi,\,\,\,e_0 := \lambda E_-,\quad e_{-0} :=
\lambda^{-1}E_+.
\eqno(2.3) 
$$ (Replacing $\hat \varphi$ by $ \varphi$ amounts to taking the quotient
by the ideal generated by the central element $c$.)  

The classical limit of $R$ involves a parameter $\hbar$. We replace
$$
 x \rightarrow \sqrt \hbar x, \quad x = e_{\pm \alpha}, H_a.\eqno(2.4)  $$
and expand
$$  R = 1 + \hbar r + o(\hbar^2).
$$  After this,    the $t_n$ are all of order $\hbar$.  Note that the
truth of this last statement is not entirely trivial.

Now recall that the Yang-Baxter relation for $R$ is equivalent to the 
recursion relation
$$   [t_n, 1 \otimes e_{-\gamma}] = (e_{-\gamma} \otimes
e^{\varphi(\gamma,\cdot)}) t_{n-1}
 - t_{n-1}   (e_{-\gamma} \otimes e^{-\varphi(\cdot,\gamma)}),\,\, n \geq
1.  \eqno(2.5)
$$  To lowest order in $\hbar$ this becomes
$$
\eqalign{ [t_1, 1 \otimes e_{-\gamma}] &= e_{-\gamma} \otimes (\varphi +
\varphi^t)(\gamma),\cr  [t_n, 1 \otimes e_{-\gamma}] &= [e_{-\gamma}
\otimes 1,t_{n-1}],\,\,\, 
  n \geq 2,\cr}\eqno(2.6)
$$  which is the same as
$$  [1 \otimes e_{-\gamma} + e_{-\gamma} \otimes 1,r - \varphi] +  [t_1,1
\otimes e_{-\gamma}  ] = 0,\,\, \gamma = 0,\cdots,l,\eqno(2.7)
$$  with $t_1$ as in in (1.5), or
$$  [1 \otimes e_{-\gamma} + e_{-\gamma} \otimes 1,\,\,r] =
\varphi(\cdot,\gamma) \wedge  e_{-\gamma}.
$$  This result is just the classical limit of the relation 
$\Delta(e_{-\gamma})R = R\Delta'(e_{-\gamma})$, which explains why it
determines $r$.

We normalize the root vectors so that the Casimir element takes the form
$$ C = \varphi + \varphi^t + \sum E_{-i} \otimes E_i + \sum E_i \otimes
E_{-i}.\eqno(2.8)
$$ Then
$$ [e_{-\gamma}, E_{-i}] = cE_{-j} {\rm ~ implies ~ that ~ }
[E_{j},e_{-\gamma}] = cE_i,
\,\, \gamma \neq 0, \eqno(2.9)
$$ 
$$ [e_{-0}, E_{-i}] = cE_{j} {\rm ~ implies ~ that ~ } [E_{-j},e_{-0}] =
cE_i, \,\,\gamma
\neq 0,\eqno(2.10)
$$  The classical r-matrix can be expressed as a formal power series in $x
= \lambda/\mu$,
$$ r = \varphi + \psi(x)^{ab}H_a \otimes H_b + \sum f_i(x) E_{-i} \otimes
E_i + \sum g_i(x) E_i
\otimes E_{-i}.\eqno(2.11)
$$  Now it is easy to work out the implications of Eq.(2.7), namely, first
taking $\gamma
\neq 0$,
$$
\eqalign{
 0 &=  [1 \otimes e_{-\gamma} + e_{-\gamma} \otimes 1, \psi(x)^{ab}H_a
\otimes H_b + \sum f_i(x) E_{-i} \otimes E_i \cr &\hskip1in + \sum g_i(x)
E_i \otimes E_{-i}] + \sum e_{-\alpha}\otimes [e_\alpha,e_{-\gamma}]\cr 
&= e_{-\gamma} \otimes \bigl(
\psi(\gamma,\cdot) + (1-f_\gamma) (\varphi + \varphi^t)(\gamma) \bigr) \cr 
&\quad \quad +
\bigl(\psi(\cdot,\gamma) - g_\gamma (\varphi + \varphi^t)(\gamma)\bigr)
\otimes e_{-\gamma}\cr &\quad \quad+
\sum f_i[e_{-\gamma},E_{-i}] \otimes E_i + {\sum}' f_iE_{-i} \otimes
[e_{-\gamma},E_i]\cr
 &\quad \quad + {\sum}' g_i[e_{-\gamma},E_i] \otimes E_{-i} + \sum g_i E_i
\otimes [e_{-\gamma},E_{-i}], \,\,\gamma
\neq 0.
\cr}\eqno(2.12)
$$  The prime on ${\sum}'$ means that the summation is over  roots that
are not simple. Cancellation  in the last two lines  imply, in view of
(2.9) and since the adjoint action is irreducible, that $f_i = f,\,g_i =
g, \, i = 1,\ldots,l$. Cancellation in the two first lines now tells us
that $\psi \propto 
\varphi + \varphi^t$, hence $\psi$ is symmetric, and it follows that $ g =
f-1$.  This gives us 
$$ r = \varphi + \sum E_{-i} \otimes E_i + g(x) C, \eqno(2.13)
$$ which is actually obvious: The two first terms is a special solution
and the last term is the only thing that commutes with
$\Delta_0(e_{-\gamma}) = 1 \otimes e_{-\gamma} + e_{-\gamma} \otimes 1$.
Next, Eq.(2.7) with $\gamma = 0$,
$$
\eqalign{
 0 &=  [1 \otimes e_{-0} + e_{-0} \otimes 1,\, \psi(x)^{ab}H_a \otimes H_b
+ \sum f_i(x) E_{-i} \otimes E_i \cr  &\hskip1in + \sum g_i(x) E_i \otimes
E_{-i}] + \sum e_{-\alpha}\otimes [e_\alpha,e_{-0}]\cr  &= E_+ \otimes
\bigl({1 \over \mu}\psi(0,\cdot) + ({1 \over \mu} - {g \over
\lambda})(\varphi +
\varphi^t)(0)\bigr)\cr &\quad \quad + \bigl( {1 \over
\lambda}\psi(\cdot,0) - {f
\over\mu}(\varphi +
\varphi^t)(0)
\bigr) \otimes E_+\cr & \quad \quad +  {f \over \mu}\sum_{i \neq +} [E_+,
E_{-i}] \otimes E_i  
  + {g \over \lambda} \sum_{i \neq +} E_{i}
\otimes [E_+,E_{-i}].\cr }\eqno(2.14)
$$  This yields $g = xf$ and the result is that
$$
 r = \varphi + \sum E_{-i} \otimes E_i + {x \over 1-x} \, C,\quad x =
\lambda/\mu,
\eqno(2.15)
$$  which agrees with the simplest r-matrix  in [BD], but in the notation
of [J].  

\vskip1in

\no {\bf 3. Twisted loop algebras.}

The construction of a twisted affine Kac-Moody algebra [K] involves two
simple Lie algebras,
${\cal L}$ and a subalgebra ${\cal L}_0$, such that ${\cal L}$ admits a
diagram automorphism of order $k = 2$ or $3$  to which is associated   a
Lie algebra automorphism $\mu$ that centralizes  ${\cal L}_0$. The
eigenvalues of $\mu$ are of the form $\omega^j,\, j = 0,1,\ldots $ , and
${\cal L} = \sum_{j=0}^{k-1} {\cal L}_{j}$, where ${\cal L}_{j}$ is the
sum of the eigenspaces  with eigenvalues $ \omega^{j\,{\rm mod}\,k}$.
 The restriction of the adjoint action of
${\cal L}$ to ${\cal L}_0$ acts irreducibly on each ${\cal L}_j$.

Now let $\{H_a, e_{\pm \alpha}\} \alpha = 1,\ldots n$ be a Chevalley basis
for 
${\cal L}_0$, and let $E_+$ be a highest weight vector (for the action of
${\cal L}_0$) in
${\cal L}_1$. Then $\{e_\alpha\},E_-$ generate ${\cal L}$, and
$$ [e_\alpha,E_+] = 0 = [e_{-\alpha},E_-].\eqno(3.1)
$$
 The twisted loop algebra $\hat {\cal L} = \Cit [\lambda, {1 \over
\lambda}] \otimes {\cal L}$ is generated by
$\{e_{\pm \alpha}\}, \alpha = 0,\ldots ,n$, with
$$  e_0 = \lambda E_-,\quad e_{-0} = {1 \over \lambda}E_+.
\eqno(3.2)$$ This algebra is of the type ${\cal A}_{cl}'$, so our standard
R-matrix applies. We define $r$ in terms of the expansion of $R$ in powers
of $\hbar$ and work out the implications of the relations (2.7).

Let $\{E_i\}$ be a Weyl basis for ${\cal L}_0$ and normalize so that the
Casimir element for that algebra is
$$ C_0 = \varphi + \varphi^t + \sum E_{-i} \otimes E_i + \sum E_i \otimes
E_{-i}.\eqno(3.3)
$$ Then a special solution of (2.7) with $\gamma \neq 0$ is given by the
first two terms in (2.13) and the general solution is
$$ r = \varphi + \sum E_{-i} \otimes E_i + \sum_0^{k-1} f_j C_j,
$$ where $C_j$ is the projection of the Casimir element $C$ of ${\cal L}$
on ${\cal L}_j$, on the first factor. Now (2.7), with $\gamma = 0$:
$$
\eqalign{ 0 &= [ 1 \otimes e_{-0} + e_{-0} \otimes 1, \, \sum E_{-i}
\otimes E_i + 
\sum f_j C_j] + \sum e_{-\alpha} \otimes [e_\alpha, e_{-0}]\cr &= { 1
\over \mu}\sum [E_+,E_{-i}] \otimes E_i + \sum f_j \bigl( {1\over \lambda}
[1 \otimes E_+, C_j] + {1 \over \mu}[E_+ \otimes 1, C_j]\bigr)\cr
&\hskip1in+ {1 \over \mu} E_+ \otimes (\varphi + \varphi^t)(0)\cr &= {1
\over \mu}[E_+ \otimes 1,C_o] + \sum \bigl({f_j \over \lambda} [ 1\otimes
E_+, C_j] + {f_{j-1}\over \mu} [ E_+ \otimes 1, C_{j-1}].
\cr}\eqno(3.4)
$$ This vanishes iff
$$
\eqalign{ & f_1 = x(f_0 +1),\,\, f_0 = x f_1,\,\,\, k = 2,\cr & f_1 =
x(f_0 +1),\,\, f_2 = x f_1, \,\, f_0 = x f_2,\,\,\, k = 3,\cr }
$$ That is,
$$ f_j = {x^j \over 1 - x^k } \, C_j - \delta_j^0 \, C_0.
$$ Finally, the unique solution is 
 $$ r = \varphi + \sum E_{-i} \otimes E_i - C_0 + {1 \over
1-x^k}\sum_0^{k-1} x^j C_j,\eqno(3.5)
$$ again in agreement with [BD], in the notation of [J].
\vskip.50cm
\no {\bf Remark.} Choose a basis of weight vectors in ${\cal L}_1$, then 
$$ C_1 = E_- \otimes E_+ + E_+ \otimes E_- + \ldots ,
$$ with unit coefficients for the contributions with highest weight. This
follows from the normalization in (3.3) and fact that $1 \otimes E_+ + 
E_+ \otimes  1 $ commutes with $C = \sum C_j$.

\vskip1in

\no {\bf 4. Including the central extension.}

\underbar {The untwisted case.} The extension is recovered by omitting the
replacement of
$\hat\varphi$ by $\varphi$ in (2.3). We can still represent the r-matrix
as a power series in $x = \lambda/\mu$, but it is no longer true, as it
was in the case of the loop group, that
$[e_0,e_{-0}] = [E_-,E_+]$. Instead,
$$ [e_0,e_{-0}] = (\hat \varphi + \hat \varphi^t)(0) = [E_-,E_+] +
c.\eqno(4.1)
$$ More generally, for polynomials $f,g \in \Cit[\lambda,{1 \over
\lambda}]$, and $x,y \in {\cal A}'_{cl}$,
$$ [fx,gx] = fg[x,y] + c\,<x,y>\, {\rm Res} (f'g),\eqno(4.2)
$$ where the form $<,>$ is the invariant form on ${\cal A}'_{cl}$
normalized as follows: If the Casimir element is $C^{ij} x_i \otimes x_j$,
then $<x_i,x_j> = (C^{-1})_{ij}$; Res$(f)$ is the constant term in
$\lambda f$.
\vskip.50cm

\no {\bf Remark.} This normalization implies that
$$ [fC_{12},gC_{23}] = fg [C_{12},C_{23}] + c_2 C_{13} {\rm Res}
(f'g).\eqno(4.3)
$$ 

\vskip.50cm

This change leaves (2.12) and (2.13) unaffected, while (2.14) becomes 
$$
\eqalign{ 0 &= E_+ \otimes \bigl({1 \over \mu}\psi(0,\cdot) + {1 \over
\mu} (\hat \varphi +
\hat \varphi^t)(0) + [e_0,g(x)E_-]\bigr)\cr  &\quad \quad + \bigl( {1
\over\lambda}\psi(\cdot,0) - {f
\over\mu}(\varphi +
\varphi^t)(0)
\bigr) \otimes E_+\cr & \quad \quad +  {f \over \mu}\sum_{i \neq +} [E_+,
E_{-i}] \otimes E_i  
  + {g \over \lambda} \sum_{i \neq +} E_{i}
\otimes [E_+,E_{-i}].\cr } 
$$  The modification in the second term ($\varphi$ replaced by $\hat
\varphi$) is exactly compensated by a new contribution from the linear
$\lambda$-term in $g$. (There  is no linear $\mu$-term in $f$.) The
conclusion is that the new r-matrix is
$$
 \hat r = \hat \varphi + \sum E_{-i} \otimes E_i + {x \over 1-x} \, C.
\eqno(4.4)
$$ 
\no \underbar {The twisted case.} It is easy to verify, with the help of
the remark at the end of Section 3, that the restitution $\varphi
\rightarrow \hat \varphi$ can be made without affecting the cancellations;
so the result is that
 $$
\hat r = \hat \varphi + \sum E_{-i} \otimes E_i - C_0 + {1 \over
1-x^k}\sum x^j C_j.\eqno(4.5)
$$

It is amusing to verify directly that the classical Yang Baxter relation
for $r$,
$$ YB(r) := [r_{12}, r_{13} + r_{23}] + [r_{13}, r_{23}] = 0,
$$ implies the same relation for $\hat r$: The inclusion of the extra term
in $\hat \varphi$ means that
$$ YB(\hat r) = YB(r) + [r_{13},(c \otimes d)_{23}].\eqno(4.6)
$$ The evaluation of $YB(r)$ now has to take into account the new term
(involving $c$) in  Eq.(4.2). Actually, only $[r_{12},r_{23}]$ is
affected, and with the aid of Eq.(4.3) one finds that the new contribution
is
$$ YB(r) = c_2 \lambda {d \over d\lambda} r_{13},
$$  which exactly cancels the other term. In the twisted case one must use
the following generalization of Eq.(4.3):
$$ [f{C_j}_{12},g{C_{j'}}_{23}] = fg [{C_j}_{12},{C_{j'}}_{23}] +
\delta_j^{j'}c_2
\,{C_j}_{13} {\rm Res} (f'g).\eqno(4.7)
$$

\vskip1in

\no {\bf 5. Deformations.}

A deformation of the standard R-matrix is a formal series
$$
 R_\epsilon=R+\epsilon R_1+\epsilon^2 R_2 + \ldots~. \eqno(5.1)
$$
\no Here $R$ is given by Eq. (1.1), with the coefficients determined by
the Yang-Baxter relation, and we determine $R_1,R_2,\ldots$ so that
$R_\epsilon$ will satisfy the same relation to each order in
$\epsilon$.  
 
Recall  that $R$ is driven by the linear term; that is, by virtue of the
Yang-Baxter relation, $R$ is completely determined by the term
 $t_1$.  It is natural to study deformations that are driven by a similar
term, with fixed but non-zero weight:
$$  R_1=S(e_{\pm\sigma} \otimes e_{\pm\rho}), \eqno(5.2)
$$
\no with $\sigma,\rho$ fixed and the factor $S$ is in ${\cal{A}'}^0$.

The following result was obtained. Among the possibilities in (5.2) the
only one that turns out to lead to a general class of deformations is 
$$
 S(e_\sigma \otimes e_{-\rho}) + S' (e_{-\rho} \otimes e_\sigma) ,\,\,\,S
\in   {\cal{A}}^0. \eqno(5.3) 
$$ 
\vskip.50cm
\no {\bf Proposition 5.1.} [F]  Suppose that $R+\epsilon R_1$ is a first
order deformation, satisfying the Yang-Baxter relation to first order in
$\epsilon$.  Suppose also that the simplest term in $R_1$   has the form
(5.3), with $S \neq 0$, then  the parameters satisfy
$$  e^{\varphi(\cdot,\rho)+\varphi(\sigma,\cdot)}=1~. \eqno(5.4)
$$  In this case  $R_1$ is uniquely determined and has the expression
$$  R_1=   (Ke_{-\rho}\otimes Ke_\sigma)R -R(Ke_\sigma\otimes Ke_{-\rho}) 
~, \eqno(5.5)
$$
\no with $K_\rho ~:=~e^{\varphi(\cdot,\rho)} = K^\sigma~:=
e^{-\varphi(\sigma,\cdot)}$.

\vskip.50cm

Deformations of this type, involving a single pair $(\rho, \sigma)$ for
which (5.4) holds, is called an elementary deformation. To first order in
$\epsilon$, the problem being then linear, one obtains a more general
space of deformations by adding the  contributions of several  such pairs,
$$ R_1 =  \sum_{(\sigma,\rho) \in [\tau]} 
\,\,\bigl(f_{-\rho}\otimes f_\sigma R -Rf_\sigma\otimes f_{-\rho} \bigr).
\eqno(5.6)
$$  Here the sum is over a subset   $[\tau]$ of the  pairs  $(\sigma,
\rho); \,\,\sigma \in \hat
\Gamma_1,\,\,\rho\in\hat\Gamma_2 $, where $\hat \Gamma_{1,2}$ are subsets
of the set of  positive generators, and
$$ e^{\varphi(\sigma,\cdot) + \varphi(\cdot,\rho)} = 1,\quad (\sigma,\rho)
\in [\tau].\eqno(5.7)
$$ But not all such compounded, first order deformations are
approximations to exact deformations  (deformations to all orders in
$\epsilon$).  

The deformed co-product was also calculated to first order in $\epsilon$,
and the result suggests an approach to the 
 exact deformations. For the following result ${\cal A}'$ is any
coboundary Hopf algebra.
\vskip.5cm

\no {\bf Theorem 5.1.} Let $R$ be the R-matrix, $ \Delta$  the coproduct,
of a coboundary Hopf algebra ${\cal A}',$ and $F \in {\cal A}'\otimes
{\cal A}'$, invertible, such that
$$
\bigl((1 \otimes \Delta_{21}) F\bigr)F_{12} = \bigl( (\Delta_{13} \otimes
1) F\bigr) F_{31}.
\eqno(5.8)
$$ Then
$$
\tilde R := (F^t)^{-1} R F\eqno(5.9)
$$ (a) satisfies the Yang-Baxter relation and (b) defines a Hopf algebra
$\tilde {\cal A}$ with the same product and with co-product
$$
\tilde \Delta = (F^t)^{-1} \Delta F^t.\eqno(5.10)
$$

\vskip.5cm
\no {\bf Proof.} (a) We substitute (5.8) into the expression $\tilde
R_{12}\tilde
  R_{13}\tilde R_{23}$. Then  use (5.7) to express $F_{12}(F_{31})^{-1}$
in terms of the co-products, and the intertwining property of $R$ ($\Delta
R = R \Delta'$) to shift the latter to the  ends. The rest is obvious. (b)
It is clear that $\tilde \Delta$ is an algebra homomorphism. We shall show
that the twisted coproduct defined by $\tilde
\Delta$ is co-associative:
$$
\eqalign{ (1 \otimes \tilde \Delta_{23}) \tilde \Delta(x) &= F^{-1}_{32}(1
\otimes \Delta_{23})
\tilde \Delta(x)F_{32}\cr &= F^{-1}_{32}(1 \otimes \Delta_{23}
\,(F^t)^{-1})(1 \otimes
\Delta_{23}\, \Delta(x))(1  \otimes \Delta_{23}\, F^t) F_{32},\cr (\tilde
\Delta_{12} +1) \tilde \Delta(x) &= F^{-1}_{21}(\Delta_{12} \otimes 1
\,(F^t)^{-1})(\Delta_{12} \otimes 1 \,\Delta(x))(\Delta_{12} \otimes 1
\,F^t) F_{21}.
\cr}
$$
   Comparing the factors at either end one gets, in view of the
co-associativity of
$\Delta$,  a relation that reduces to (5.8) after  re-numbering the
spaces. The theorem is proved. \footnote*{The connection between Eq.(5.8)
and co-associativity was pointed out to  me by Kajiwara.}    
\vskip.5cm

The naturality of this construction is indicated by the following simple
result.
\vskip.5cm
\no {\bf Proposition 5.2.} Let $R,F$ be as in Theorem 5, and suppose that
there is a second twist $\tilde F
\in \tilde {\cal A} \otimes \tilde {\cal A}$, invertible, such that 
$$
\bigl((1 \otimes \tilde \Delta_{21}) \tilde F\bigr)\tilde F_{12} = \bigl(
(\tilde \Delta_{13} \otimes 1) \tilde F\bigr) \tilde F_{31}.
$$ Then the two twists can be combined to $G = F\tilde F$ satisfying
$$
\bigl((1 \otimes \Delta_{21}) G\bigr)G_{12} = \bigl( (\Delta_{13} \otimes
1) G\bigr) G_{31}.
$$

 We return to our subject, with $R$ again denoting the standard R-matrix
of the algebra
${\cal A}' = {\cal A}/I$. We show first that interesting solutions of
(5.8) exist.  Then we do some preliminary calculations that help us make a
general ansatz for $F$ in the form of a double expansion, $ F = \sum 
\epsilon^{nm}F_n^m$, and finally we derive a recursion relation for
$F_n^m$ that will allow us to calculate the classical limit. 
\vskip.5cm
\no {\bf Examples.} An exact deformation of $R$, with first order term
$R_1$ as in (5.5) is given by
$$ F =  e_q^{-\epsilon f_\sigma \otimes f_{-\rho}},\eqno(5.11) 
$$ with
$$ f_\sigma := e^{-\varphi(\sigma,\cdot)}e_\sigma,\,\,\, f_{-\rho} :=
 e_{-\rho}e^{\varphi(\cdot,\rho)}\eqno(5.12)
$$ The q-exponential is   as follows: 
 $q = e^{\varphi(\sigma,\rho)}$, ${e_q}^A := \sum A^n/[n!]_q,\,\, [n!]_q =
[1!]_q \ldots [n]_q,\,\, \break [n]_q = (q^n-1)/(q-1)$.  Note that, if $AB
= qBA$, then $e_q^A e_q^B = e_q^{(A+B)}$. This construction works easily
for some compound deformations, with (5.5) replaced by (5.6). Proposition
5.2 shows that an elementary twist $F$, of the simple form (5.11), can be
combined with another elementary twist $\tilde F$, of the same  type but
with $(\sigma,\rho)$ replaced by $(\sigma',\rho')$, only if
$\tilde\Delta(f_\sigma'), \tilde\Delta (\rho')$ reduce to
$\Delta(f_\sigma'),
 \Delta (\rho')$; that is, only when the four generators quommute among
themselves. 

\vskip.5cm
\no {\bf Notation.} From now on it will be convenient to use the
generators $f_{\pm \alpha}$ defined in (5.12). The standard co-product
then takes the form
$$
\Delta f_\sigma =  K^\sigma \otimes f_\sigma + f_\sigma \otimes 1, \quad
\Delta f_{-\rho} = 1 \otimes f_{-\rho} + f_{-\rho} \otimes K_\rho,
$$ with
$$ K_\rho := e^{\varphi(\cdot,\rho)},\,\,\,K^\sigma :=
e^{-\varphi(\sigma,\cdot)}.
$$

The general case of compound deformations is much more complicated. The
calculations are manageable only so long as $F$ can be constructed from 
elements of the type $f_\sigma
\otimes f_{-\rho}$ only,  with the factors in this order. A general result
is  Theorem 5.2 below. We need some preparation.
 
\vskip.5cm
 \no {\bf Proposition 5.3.} Let $R_\epsilon$ be an exact deformation of
the type   
$$
\eqalign{ &R_\epsilon = (F^t)^{-1}RF, \quad \quad F =  
\sum \epsilon^n  (F_n + \ldots ),\cr &F_n = \sum_{(\sigma,\rho) \in
[\tau]} F_{(\sigma)}^{(\rho')} f_{\sigma_1}\ldots f_{\sigma_n}  \otimes
f_{-\rho_1'} \ldots f_{-\rho_n'},\cr}
\eqno(5.13)
$$ where $+ \ldots $ stands for terms with less than $n$ factors.
  Let
$\Gamma_1, \Gamma_2$ be the subalgebras of ${\cal A}'^+$ generated by
$\hat\Gamma_1,
\hat\Gamma_2$. Then we have:  (a)  There is an isomorphism $\tau: \Gamma_1
\rightarrow \Gamma_2$, such that the set $[\tau]$ is the restriction of
the graph of $\tau$ to $\hat\Gamma_1,
\hat\Gamma_2$,
$$ [\tau] = \{\sigma, \rho \,|\, \sigma \in \hat\Gamma_1, \, \rho = \tau
\sigma\in
\hat\Gamma_2\}.\eqno(5.14)
$$ (b) The elements $F_n$ satisfy the recursion relations
$$
 [F_n, f_{-\sigma} \otimes 1] =  ( K^\sigma \otimes f_{-\rho})F_{n-1}
 - F_{n-1}(K_\sigma
\otimes f_{-\rho}),  \quad  (\sigma,\rho) \in [\tau],\eqno(5.15)
$$
\vskip-2mm
\no  as well as
\vskip-4mm
$$
 [1 \otimes f_\rho, F_n] = F_{n-1}\bigl( f_\sigma \otimes K^\rho
\bigr) - \bigl(f_\sigma \otimes K_\rho\bigr) F_{n-1} .\eqno(5.16)
$$
 (c) These recursion relations have the unique solution  
$$ F_{(\sigma)}^{(\rho')} = -\overline t_{(\sigma)}^{(\sigma')},\,\,\,\,
(\rho'_1, \ldots ,
\rho'_n) =
\tau(\sigma'_1, \ldots , \sigma'_n),\eqno(5.17)
$$ where the coefficients on the right are the same as in Eq.(1.1), except
that $\varphi$ is replaced by $-\varphi: \,\,t(\varphi) = \overline
t(-\varphi)$.  
\vskip.5cm 
\no {\bf Proof .} We begin by offering some justification for the
assumptions.  In view of the form of $R_1$ it is expected that
$R_n$ is a sum of products of factors of three types:
$$ e_{-\alpha} \otimes e_\alpha,\,\,\,f_{-\rho} \otimes f_{\sigma},\,\,\, 
f_{\sigma} \otimes f_{-\rho},\quad \sigma \in \hat\Gamma_1,\, \rho \in 
\hat\Gamma_2,\eqno(5.18)
$$  with coefficients in ${\cal A}' \otimes {\cal A}'$.  In $R_n$ , we
isolate the terms with the highest number of factors of the third type, 
$$ X_n = \sum A_{(\sigma)}^{(\rho')}\bigl(e_{-\alpha_1} \ldots
e_{-\alpha_k} \otimes e_{\alpha_1} \ldots e_{\alpha_k}\bigr)
B_{(\sigma)}^{(\rho')} \bigl(f_{\sigma_1} \ldots f_{\sigma_n} \otimes
f_{-\rho_1'} \ldots f_{-\rho_n'}\bigr),  
$$
   We shall show that $R_n$contains $X_n \neq 0$.

Let 
$$ Y\hskip-1mm B_\epsilon := R_{\epsilon \,12} R_{\epsilon\,13}
R_{\epsilon\,23} -  R_{\epsilon\,23} R_{\epsilon\,13} R_{\epsilon\,12} \in
{\cal A}' \otimes {\cal A}' 
\otimes {\cal A}'. \eqno(5.19)
$$ All terms in $Y\hskip-1mmB_\epsilon$ of order $\epsilon^n$, that have
$n$ factors of the third type in spaces 1,2 are contained in
$$ P_n := F_{n\,12} R_{13} R_{23} - R_{23} R_{13} F_{n\,12}.\eqno(5.20)
$$ For these terms to cancel among themselves $X_n$ must take the form 
$$ X_n = R F_n,\quad F_n = F_{(\sigma)}^{(\rho')} \bigl(f_{\sigma_1} \ldots
f_{\sigma_n} \otimes f_{-\rho_1'} \ldots f_{-\rho_n'}\bigr).\eqno(5.21)
$$ The sum is over all pairs  $(\sigma,\rho) \in [\tau]$ and all
permutations $(\rho')$ of
$(\rho)$.
 
Next, the recursion relation (5.15) follows easily from the Yang-Baxter
relation ( more precisely from an examination of terms of low order in
space 2), and (5.16) from a similar calculation. (Recursion relations of
this type were examined in detail in [F], so we skip the details.)
 
We have $F_0 = 1$ and $F_1 = \sum f_\sigma \otimes f_{-\rho}$. Taking
$n=1$ in (5.15) or (5.16) one gets,  
$$ [f_\alpha, f_{-\beta}] =  \delta^\beta_\alpha
\bigl(e^{\varphi(\cdot,\alpha)} -
e^{-\varphi(\alpha,\cdot)}\bigr),\eqno(5.23)
$$ which is confirmed by the definitions in (5.12) and the relation (1.4).
When (5.15) is reduced to a recursion relation for the coefficients, then
it turns out to agree (up to the sign of $\varphi$) with the recursion
relations for the coefficients 
$t_{(\sigma)}^{(\sigma')}$ that is implied by (1.1) and (2.5).  The
integrability of these relations is precisely the statement (a) of  the
theorem, as follows easily from the analysis of these recursion relations
in [F]. Finally,  when (a) holds, then the  relation (5.16) is equivalent
to (5.15). The proposition  is proved.
\vskip.5cm After these preliminary explorations we are able to formulate a
general result.   
\vskip.5cm
\no {\bf Theorem 5.2.} Let $\Gamma_1, \Gamma_2$ be subalgebras of ${\cal
A}'^+$, generated by subsets $\hat \Gamma_1,\hat \Gamma_2$ of the
generators, and $\tau: 
\Gamma_1 \rightarrow \Gamma_2$ an algebra isomorphism.   Let $F \in {\cal
A}' \otimes {\cal A}'$ be a formal series of the form
$$ F = 1 + \sum_{n = 1}^\infty \epsilon^n F_n, \quad F_n = \sum
F_{(\sigma)}^{(\rho)}  f_{\sigma_1}\ldots f_{\sigma_n} \otimes f_{-\rho_1}
\ldots f_{-\rho_n}.
\eqno(5.24)
$$ The second sum is here over all $\sigma_i \in \hat \Gamma_1,  \rho_i
\in \hat\Gamma_2$. (!) Note that $F_n$ is   a power series in $\epsilon$.
Suppose that $F$ satisfies (5.8), then
$$ F_1 = -\sum_{\tau^m \sigma = \rho}\epsilon^{m-1}(f_\sigma \otimes
f_{-\rho}),\eqno(5.25)
$$ and
$$ (1 \otimes K_\rho \partial_\rho) F_n + \sum_{\tau^m \sigma = \rho} 
\epsilon^m \,\, [1 \otimes f_\sigma,F_n] + \sum_{\tau^m \sigma =
\rho}\epsilon^{m-1}(f_\sigma
\otimes K^\sigma) F_{n-1} = 0.\eqno(5.26)
$$ This last relation is satisfied for $n = 1$ by (5.25), with $F_0 = 1$.
With $F_1$ thus fixed,  $F_2,F_3,\ldots$ are determined recursively and
uniquely.
\vskip.5cm
\no {\bf Notation.} The sums in (5.25-6), and similar sums to follow,
should be understood  to run over $\sigma \in \hat \Gamma_1$ and over all
values of the integer $m$ such that $\tau^m\sigma$ is defined; that is,
all values of $m$ such that $\tau^{m-1}\sigma \in \hat\Gamma_1$.
\vskip.5cm
\no {\bf Proof.} That Eq.(5.8) implies (5.25) and (5.26) is a simple
calculation; one collects all terms that have exactly one generator in the
second space.   Let  us verify that the recursion relation is satisfied for
$ n = 1$ by (5.25). The second term is
$$ -\sum_{\tau^{m'} \sigma'  = \rho'}\epsilon^{m'} 
 \sum_{\tau^m \sigma = \rho}\epsilon^{m-1}f_{\sigma'} \otimes [f_{\sigma},
f_{-\rho'}].
$$ The commutator is
$$
 [f_{\sigma}, f_{-\rho'}] = e^{\varphi(\cdot,\rho')} -
e^{-\varphi(\rho',\cdot)}
 = e^{\varphi(\cdot,\tau^{m'}\sigma')} -
e^{\varphi(\cdot,\tau^{m'+1}\sigma')}.
$$ The double sum reduces to $\sum_{\tau^m \sigma =
\rho}\epsilon^{m-1}f_\sigma
\otimes (K^\sigma - K_\rho)$ and (5.26) reduces to an identity. It remains
to prove that  (5.26) has a unique solution.  Consider first the case that
$\hat\Gamma_1 \cap \hat\Gamma_2
$ is empty; then the second term in (5.26) vanishes and the third term
reduces to the term $m=1$. The recursion relation then reduces to the same
form as that which determines the coefficients of the standard R-matrix,
which is known to be integrable [F]. (In this case Proposition 5.3 is the
complete solution of the problem, for there are no terms ``$+\ldots$" in
(5.13).) In the general case, when 
$\hat\Gamma_1 \cap \hat\Gamma_2
$ can be non-empty, the second term in (5.26) makes the solution  more
difficult,  but the existence of a solution can still be proved. To do
this we  expand $F_n$ as a power series in $\epsilon$, with constant term

$$ F_n^1 = \sum F_{(\sigma)}^{(\tau\sigma')} f_{\sigma_1'}\ldots
f_{\sigma_n'} \otimes  f_{-\tau\sigma_1'}\ldots f_{-\tau \sigma_n'},
$$ and determine the coefficients recursively. The problem is therefore 
always the integrability of $K_\rho\partial_\rho X = Y, \rho \in
\hat\Gamma_2$, with $Y \in {\cal A}'$  
 given, and this is known [F] to have a unique solution in ${\cal A}'$, as
already noted.  The theorem is proved.  
\vskip.5cm The converse, that the solution of (5.26) with $F_0=1$ and
$F_1$ given by (5.25) satisfies (5.8) (and therefore gives a solution of
the Yang-Baxter relation) was proved only in the special case that 
$\hat\Gamma_1 \cap \hat\Gamma_2$ is empty.   Note, however, that the exact
form (5.25) of
$F_1$  can be inferred directly from the Yang-Baxter relation. Further
direct computation  supports the idea that $R_\epsilon$ always has the
form $(F^t)^{-1}RF^t$, with F of the  form assumed in (5.24). This is
strong support for the belief that the solution of the  recursion relation
(5.26), which was proved to exist always, actually furnishes the solution
to the problem of exact deformations in the general case. As we shall see,
additional favorable evidence comes from an examination of the classical
limit. To prepare for this we need
\vskip.5cm
\no {\bf Proposition 5.4.} Let 
$$ F_n^m = \sum_{\rho = \tau^m \sigma}  \overline t_{(\sigma)}^{(\sigma')}
f_{\sigma_1}
\ldots f_{\sigma_n}
\otimes f_{-\tau^m \sigma_1'} \ldots f_{-\tau^m \sigma_n'},\quad F^m_0 =
1, 
 \eqno(5.27)
$$ in which the sum extends over $\sigma_i \in \hat \Gamma_1$, $(\sigma')$
a permutation of $(\sigma)$, and the coefficients $\overline
t_{(\sigma)}^{(\sigma')}$ are the same as in (5.17). Then the unique
solution of (5.26) is
$$
\eqalign{ F_n &= \sum_{\Sigma n_i = n} \epsilon^{n_2 + 2n_3 + \ldots}
F_{n_1}^1F_{n_2}^2\ldots\ =  F_n^1 - \epsilon F_{n-1}^1 F^2_1 +
\epsilon^2\bigl(F^1_{n-2} F^2_2 + F_{n-1}F^3_1\bigr)
 + \ldots,\cr F &= \sum \epsilon^nF_n = \sum \epsilon^{n_1 + 2n_2 +
\ldots}F^1_{n_1}F^2_{n_2}
\ldots = F^1F^2\ldots,\quad F^m = \sum \epsilon^{nm} F_n^m.\cr}\eqno(5.28)
$$
\ve 
\no {\bf 6. Esoteric r-matrices.}

 We specialize to the case of a quantized Kac-Moody algebra of finite
type.  
\vskip.5cm
\no {\bf Proposition 6.1.} If ${\cal A}'$ is a quantized Kac-algebra of
finite type, then $\hat\Gamma_1$ is a proper subset of the set of positive
generators and $\tau^{m+1}\hat\Gamma_1 \cap \hat\Gamma_1$ is a  proper
subset of $\tau^m\hat\Gamma_1 \cap \hat\Gamma_1$.
\vskip.5cm
\no {\bf Proof.} Suppose that the statement is false. Then there is
$f_\sigma \in \hat\Gamma_1$ such that 
$\tau^m f_\sigma \in \hat\Gamma_1$ for all $m$, and consequently $
\tau^kf_\sigma = f_\sigma$ for some k. But the condition (5.7), in the
classical limit, implies that
$$
\varphi(\tau^m\sigma,\cdot) + \varphi(\cdot.\tau^{m+1}\sigma) = 0.
$$ Summing over $m=0,1,\ldots,k-1$ we obtain
$$
\sum_m(\varphi + \varphi^t)(\tau^m\sigma) = 0,
$$ which contradicts the fact that the Killing form is non-degenerate.
\vskip.5cm In the classical limit
$$ R_\epsilon = 1 + \hbar r_\epsilon + o(\hbar^2)\quad r_\epsilon = r +
\epsilon + o(\epsilon^2).\eqno(6.1)
$$ In the case of an exact elementary deformation $R_\epsilon$,   the
associated  exact deformation $r_\epsilon$ of $r$ coincides with the first
order,
$$ r_\epsilon = r + \epsilon r_1.  \eqno(6.2)
$$

Consider the general case of an exact deformation of $R$ of the form
postulated in   Theorem 5.2. Define $X_\epsilon$  by  
$$ F = 1 + \hbar X_\epsilon + o(\hbar^2), \eqno(6.3)
$$ so that 
$$ r_\epsilon = r + X_\epsilon - X^t_{\epsilon}.\eqno(6.4)
$$
\vskip.5cm
\no {\bf Notation.} In this section the symbols $\Gamma_{1,2}$ stand for
Lie algebras,  the classical limits of the algebras so designated until
now.
\vskip.5cm

From the fact that the coefficients in (5.27) are the same as the
coefficients in (1.1), and the known classical limit of the standard
R-matrix for a Kac-Moody algebra of finite type, we get without
calculation that
$$ X_\epsilon = -\sum_m \sum_{E_i \in \Gamma_1}\epsilon^{nm} E_i \otimes
E_{-\tau^m i},\eqno(6.5)
$$ in which $n$ is the height of $E_i$. The normalization is the same as in
Sections 2-4; more precisely it is fixed as follows. (a) The set $\{E_i\}$
includes the  generators of $\Gamma_1$. (b) The statement (2.9).
 \footnote*{ Condition (b) can be re-phrased as follows. Let $\Gamma_1^-$
be the  Lie algebra generated by $\{f_{-\sigma}\},~f_\sigma \in \hat
\Gamma_1$ and $\Gamma$ the  Lie  algebra generated by $\{f_{\pm
\sigma}\},~f_\sigma \in \hat\Gamma_1$. Then
$$
\sum_{E_i \in \Gamma_1} E_i \otimes E_{-i}
$$ is the projection on $\Gamma_1 \otimes \Gamma^-$ of a
$\Gamma$-invariant element of 
$\Gamma \otimes \Gamma$.}  Consequently, 
$$ r_\epsilon = r -  \sum_m \sum_{\matrix {E_i \in \Gamma_1\cr E_j =
\tau^m E_i
\cr}} \epsilon^{nm}
\,\, E_i 
\wedge E_{-j}.\eqno(6.6)
$$ The sums are finite, by Proposition 6.1.   A renormalization exists
that reduces the numerical coefficients to  unity ($\epsilon$ now
interpreted as in ~\Cit); the result is in complete agreement with [BD].

\vskip.5cm

\no \underbar {Deformations in the affine case.} Let ${\cal A}'$ be a
quantized Kac-Moody algebra of affine type. Two cases should be
distinguished. If the subsets $\hat\Gamma_{1,2}$ of positive roots do no
include the  imaginary root $e_0$, then the formula (6.6) applies without
change, except that now $r$  is   one of the standard affine r-matrices
determined earlier, Eq.s~(2.15), (3.5), (4.4) or (4.5). There is nothing
more to be said about this case and we turn our attention to the other one.

 What  merits special attention is the possibility that the first order
deformation (5.6)  may include one of the following
$$ e_0 \wedge e_{-\rho} =  \mu(E_- \otimes e_{-\rho}) -\lambda (e_{-\rho}
\otimes E_-),
\eqno(6.7)
$$  or
$$ 
 e_\sigma \wedge e_{-0} = \lambda^{-1} (e_\sigma \otimes E_+) -
\mu^{-1}(E_+ \otimes e_\sigma),\eqno(6.8)
$$ with 
$$
\varphi(\cdot,\rho) + \varphi(0,\cdot) = 0, \quad {\rm resp.} \quad
\varphi(\cdot,0) +
\varphi(\sigma,\cdot) = 0,\eqno(6.9)
$$ which implies that $\rho  \neq 0$, resp. $ \sigma \neq 0$. A simple
renormalization, that connects the principal picture to the homogeneous
picture, brings   (6.8) to the form
$$ e_\sigma \wedge e_{-0} = \sqrt  {\mu/\lambda} (e_\sigma \otimes E_+) -
\sqrt{\lambda/\mu}(E_+ \otimes e_\sigma).
$$ 
  
To deal with the general case of exact deformations it is useful to note
the following
\vskip.5cm
\no {\bf Proposition 6.2} If ${\cal A}'$ is a quantized Kac-Moody of
affine type, then {\it either }the  statement about $\hat\Gamma_1$ in
Proposition 6.1 continues to hold, {\it or } ${\cal A}'$ is of type 
$A_{N-1}^{(1)}$, $\hat\Gamma_1$ consists of all the positive generators,
and $\tau$ generates the cyclic  group of order $N$.
\vskip.5cm
\no{\bf Proof.} Suppose there is $f_\sigma \in \hat\Gamma_1$ such that
$\tau^Nf_\sigma = f_\sigma $ for some $N$.  Then the Killing form is
degenerate. But it is known [K] that any subalgebra of a Kac-Moody algebra
of affine type, obtained by removing one generator, is a Kac-Moody algebra
of finite type. It follows that  
  $\hat\Gamma_1$ contains all the positive generators and exactly one
$\tau$ orbit.  Then $\hat\Gamma_1 = \hat\Gamma_2$ and 
$\tau$ lifts to an isomorphism of the Dynkin diagram, which implies the
result.
\vskip.5cm  
\no In this section we exclude the exceptional case. This means that the
classical limit of $\Gamma_1$ is a finite dimensional Lie algebra, so that
(6.6) can be applied directly, since the sum is finite.

Alternatively, the classical limit of an exact deformation can be found
with the  help of the recursion relation  
 $$ (1 + K_\rho\partial_\rho)F^m_n = -  (f_\sigma \otimes
K_\rho)F^m_{n-1},\quad 
\tau^m \sigma = \rho,\eqno(6.10)
$$ or better, the equivalent relation
$$ [1 \otimes f_\rho, F^m_n] =    -  \biggl((f_\sigma \otimes K_\rho)
F^m_{n-1} -  F^m_{n-1}(f_\sigma
\otimes K^\rho)\biggr)\eqno(6.11)
$$ for $F^m_n = \delta_n^0 + \hbar X^m_n + o(\hbar^2)$. This implies that
$X^m = 
\sum_{n = 0,1,\ldots} \epsilon^{mn} X^m_n$ (a finite sum) is the unique
solution (of the form that appears in (6.5)) of
$$ [1 \otimes f_\rho +  \epsilon^m f_\sigma \otimes 1, X^m] =  \epsilon^m  
f_\sigma \otimes  (\varphi + \varphi^t)(\rho),   \quad \tau^m\sigma = \rho
\in
\hat\Gamma_2.\eqno(6.12)
$$
 \vskip.5cm
\no {\bf Example.} Let ${\cal A}'_{cl}$ be the untwisted, affine Kac-Moody
algebra 
$\tilde {\cal L}, {\cal L} = {sl(N)}$. A set of positive Serre generators
is provided by the unit matrices $e_i = e_{i,i+1}, ~ i = 1, \ldots N-1$.
Set $e_N = e_0 = \lambda e_{N1}$. The ``most esoteric" deformation (the
one with the largest $\Gamma_1$) is defined  as follows. Take $\Gamma_1$
to be generated by $e_i,~i = 1, \ldots N-1$,  and $\tau e_i = e_{i+1},~i =
1,\ldots N-1$.\break Then $X^m = \sum_n \epsilon^{nm}X_n^m$ with 
$$ X^m_n = -\sum_{i+m+n \leq N}e_{i,i+n} \otimes e_{i+m+n,i+m}
-\sum_{i+m+n = N+1} e_{i,i+n} \otimes \lambda^{-1} e_{1,i+m} 
$$  and
$$ r_\epsilon = r + \biggl( \sum \epsilon^{nm} X^m_n -  {\rm
transpose}\biggr).
$$ Taking $N=3$ one obtains
$$ r_\epsilon = r - \biggl(\epsilon \,e_{12} \otimes e_{32} + 
\epsilon^2 e_{13} \otimes \lambda^{-1} e_{12} + \epsilon^2 e_{12} \otimes
\lambda^{-1} e_{13} - {\rm transpose}\biggr),
$$ and the renormalization $e_{ij} \rightarrow \lambda^{{j-i \over
3}}e_{ij}$ gives the final  result
$$
\eqalign{ r_\epsilon = r &- \epsilon \{ \xi^{-1}e_{12} \otimes e_{32}  +
\xi^{-1} e_{23} \otimes e_{13} + \xi^{-2} e_{13} \otimes e_{12}\}   -
\epsilon^2
\xi^{-1}e_{12} \otimes e_{13}\cr & + \epsilon \{ \xi e_{32} \otimes e_{12} 
+ \xi  e_{13} \otimes e_{23} + \xi^2 e_{12} \otimes e_{13}\} + \epsilon^2
\xi e_{13} \otimes e_{12},\cr }
$$ with $\xi = (\lambda/\mu)^{1/3}$. The un-deformed piece is
$$ r = \phi + \sum_{i<j}e_{ij} \otimes e_{ji} = {1\over 3}\big(\sum e_{ii}
\otimes e_{ii} -  e_{11} \otimes e_{22} - e_{22} \otimes e_{33} - e_{33}
\otimes e_{11}\bigr)+ \sum_{i<j}e_{ij} \otimes e_{ji},
$$
$\varphi$ being completely fixed by the relations (5.7).  This is in
agreement with [BD], after transposition and setting $\xi = e^{u/3},
\epsilon = 1$. .

\ve

\no {\bf 7. Elliptic r-matrices.} 

Here we consider the exceptional case (Proposition 6.1) in which
$\hat\Gamma_1$ contains all the generators of ${\cal A}'^+$, ${\cal A}'$
is of type $A_{N-1}^{(1)}$ and $\tau^N = 1$.

The expression (5.25) for $F_1$ can be justified as before and the sum is
convergent if we interpret $\epsilon$ in \Cit and stipulate that
$$ |\epsilon| < 1,
$$ namely
$$ F_1 = {-1\over 1 - \epsilon^N}\sum_{m=1}^N \sum_{\sigma \in
\hat\Gamma_1} \epsilon^m  f_\sigma \otimes f_{-\tau^m\sigma}.\eqno(7.1)
$$ Most, but not all, of the infinite sums that arise can be made
meaningful in this way. In particular, (5.25) becomes 
$$ (1-\epsilon^N)(1 \otimes K_\rho\partial_\rho)F_n + \sum_{m=1}^N
\epsilon^m\, [1 \otimes f_{\tau^{-m}\rho},F_n] + \sum_{m=1}^N
\epsilon^{m-1}(f_{\tau^{-m}\rho} \otimes
K^{\tau^{-m}\rho})F_{n-1}.\eqno(7.2)
$$ We verify directly that it holds for $n=1$. The second term is
$$\eqalign{ {-1\over 1-\epsilon^N}\sum_{n=1}^N
\sum_{m=1}^N&\epsilon^{m+n}f_{\tau^{-m-n}\rho}
\otimes (K_{\tau^{-m}\rho} - K_{\tau^{1-m}\rho})\, \cr &={-1\over
1-\epsilon^N}\sum_{M=1}^N \epsilon^M f_{\tau^{-M}\rho}
\otimes (K^{\tau^{-M}\rho} - K_\rho)(1-\epsilon^N).
\cr}
$$ The last factor comes from the fact that the equation $\tau^{m+n} =
\tau^M,\, \,n,M $ given, only determines 
$m$ Mod $N$. The term $K^{\sigma'}(K_\rho)$ comes from the ends of the 
summation   while all the other terms cancel pairwise since $K^\sigma =
K_{\tau\sigma}$.

The infinite product
$$ F =  F^1F^2\ldots\eqno(7.3)
$$ cannot be given anything more than a formal significance in the
structural context but, as will be shown below, in a finite  dimensional
representation the question of convergence (with $\epsilon $ in \Cit) is
not  difficult. We define $F^m$ by the (always uniquely integrable)
relation (6.10),
$$ (1 \otimes K_\rho\partial_\rho)F^m =
-\epsilon^m(f_{\tau^{-m}\rho}\otimes K_\rho)F^m,\quad F^m = 1 - \epsilon^m
\sum_\sigma f_\sigma \otimes f_{-\tau^m\sigma} +
o(\epsilon^{2m}),\eqno(7.4)
$$ or its equivalent
$$ [1 \otimes f_\sigma,F^m] = -\epsilon^m\biggl((f_{\tau^{-m}\sigma}
\otimes K_\sigma)F^m -F^m(f_{\tau^{-m}\sigma} \otimes
K^\sigma)\biggr),\eqno(7.5)
 $$   with the same initial condition. We verify that, with this
definition of $F^m$, (7.3) satisfies (7.2) or
$$ (1 - \epsilon^N)(1 \otimes K_\rho\partial_\rho)F +
\sum_{\tau^n\sigma=\rho} \epsilon^n[1 \otimes f_\sigma,F]
+\sum_{\tau^n\sigma=\rho} \epsilon^n(f_\sigma \otimes K^\sigma)F =
0.\eqno(7.6)
$$ The range of the summation is $ n =  1,2,\ldots,N,\, \sigma \in
\hat\Gamma_1$.  One has
$$
\eqalign{ &\sum_n \epsilon^n[1 \otimes f_{-\tau^{-n}\rho},F^mF^{m+1}] =
-\sum_n \epsilon^{m+n} f_{\tau^{-m-n}\rho} \otimes
K_{\tau^{-n}\rho}F^mF^{m+1}\cr &\hskip.3cm - F^m\biggl\{ - \sum_n
\epsilon^{m+n} f_{\tau^{-m-n}\rho} \otimes K^{\tau^{-n}\rho} + 
\sum_n \epsilon ^{m+n+1} f_{\tau^{-m-n-1}\rho} \otimes
K^{\tau^{-n-1}\rho}\biggr\} F^{m+1} + \ldots \,.\cr}
$$ In the second line everything cancels except for the first and the last
terms, leaving
$$ -\sum_n \epsilon^{m+n}( f_{\tau^{-m-n}\rho} \otimes K_{\tau^{-n}\rho})
F^m F^{m+1}  + (1-\epsilon^N)F^m \epsilon^{m+1}( f_{\tau^{-m-1}\rho}
\otimes K_\rho) F^{m+1} + \ldots\,.
$$ The total contribution of the commutator in (7.6) is thus
$$ -\sum_{n=1}^N \epsilon^{n+1}( f_{\tau^{-n-1}\rho} \otimes
K_{\tau^{-n}\rho})F  + (1-\epsilon^N)\sum_{m=1}^\infty F^1\ldots F^m
\epsilon ^{m+1}(f_{\tau^{-m-1}\rho} \otimes K_\rho) F^{m+1} \ldots \,.
$$ Adding the first term in (7.6) leaves us with
$$ -\sum_n \epsilon^{n+1}( f_{-\tau^{-n-1}\rho}\otimes K_{\tau^{-n}\rho})F
- 
\epsilon(1-\epsilon^N) (f_{-\tau^{-1}\rho} \otimes K_\rho) F  = -\sum_n
\epsilon^n (f_{-\tau^{-n}\rho} \otimes K^{\tau^n\rho}) F,
$$
 which is cancelled by the last term.

In the classical limit $F^m = 1 +  \hbar X^m + o(\hbar^2)$ and $X^m$
satisfies (6.12). We shall solve these relations in the case of the
simplest affine Kac-Moody algebra. Set

$$[f_1,f_{-1}] = (\varphi + \varphi^t)(1) = \sigma_3,
$$ and
$$X^m = A^m \sigma_3 \otimes \sigma_3 + B^m(f_1 \otimes f_{-1} + f_0
\otimes f_{-0}) + C^m(f_1 \otimes f_{-0} + f_0 \otimes f_{-1})
$$ and impose (6.12). The result is, with $x =  \sqrt{\lambda/\mu}$,
$$\eqalign{ A^m &= \sum_{n=1}^\infty (-\epsilon^{2n})^m x^{-n},\cr B^{2m}
&= \sum_{n=1}^\infty(\epsilon^{2n-1})^{2m}x^{1-n}, \quad B^{2m-1} = 0,\cr
C^{2m-1} &= \sum_{n=1}^\infty(\epsilon^{2n-1})^{2m-1}x^{1-n}, \quad C^{2m}
= 0.\cr}
$$ The deformed r-matrix is $r_\epsilon = r + X - X^t$, with 
$$\eqalign{ X& = \sum_{n=1}^\infty X^m = \sum_n {-\epsilon^{2n}\over 1 +
\epsilon^{2n}}x^{-n} \sigma_3\otimes \sigma_3 \cr &+ \sum_{n=1}^\infty
{\epsilon^{4n-2}\over 1- \epsilon^{4n-2}}x^{1-n}(f_1 \otimes f_{-1} + f_0
\otimes f_{-0})  + \sum_{n=1}^\infty {\epsilon^{2n-1}\over 1-
\epsilon^{4n-2}}x^{1-n}(f_1 \otimes f_{-0} + f_0 \otimes f_{-1}).\cr}
$$ Setting $\lambda/\mu = e^{2\pi iu}$ one gets
$$\eqalign{ (i/2)(X-X^t) = &\sum_{n=1}^\infty\biggl\{ {-\epsilon^{2n}\over
1 + \epsilon^{2n}}(\sigma_3 \otimes \sigma_3) \sin 2n\pi u \cr &\hskip.2cm
+ {\epsilon^{4n-2}\over 1 - \epsilon^{4n-2}}\bigl(  x \,f_1 \otimes f_{-1}
+ {1\over x}\, f_{-1} \otimes f_1\bigr) \sin (2n-1)\pi u\cr  &\hskip.5cm +
{\epsilon^{2n-1}\over 1 - \epsilon^{4n-2}}\bigl(\sqrt{1/\mu\lambda}\,f_1
\otimes f_1 + \sqrt{\mu\lambda}\,f_{-1}
\otimes f_{-1}\bigr) \sin (2n-1) \pi u\biggr\}.\cr}
$$ The trigonometric r-matrix (2.15) is
$$ {i\over 2}\biggl({1\over\tan \pi u}(\sigma_3 \otimes \sigma_3) +
{1\over   \sin \pi u} \bigl(\sqrt x\,f_1 \otimes f_{-1} +
\sqrt{1/x}\,f_{-1} \otimes f_1\bigr)\biggr).
$$ Adding, one finds the series expansion of elliptic functions, and
complete agreement with the elliptic  r-matrices of [BD]. To transform to
their notation replace 
$$f_1 \rightarrow \sqrt\lambda\, e_{12}, \,f_{-1} \rightarrow
\sqrt{1/\lambda}\, e_{21}\eqno(7.7)
$$

Finally, we shall show that the expression for the Universal Elliptic
R-matrix as an infinite product is both  meaningful and useable, by
projecting on a finite dimensional represention. We limit ourselves to
the  fundamental representation of $sl(2)$. After  rescaling of the
generators as in (7.7), $F^m$ and $R_\epsilon$ take the form
$$ F^m = \pmatrix{a^m&&&d^m\cr&b^m&c^m&&\cr &c^m&b^m&&\cr a^m&&&d^m},\quad
R_\epsilon =  \pmatrix{a&&&d\cr&b&c&&\cr &c&b&&\cr a&&&d \cr }.
$$  The matrix elements are completely determined by the recursion
relation (7.5); namely for  $ m = 1,2,\ldots $,
$$
\eqalign{ a^{2m-1} &= 1 - \epsilon^{4m-2},~~  b^{2m-1} = 1 -
\epsilon^{4m-2}{q^2\over x},
 \,\, c^{2m-1} = 0, ~~d^{2m-1} = \epsilon^{2m-1}({1\over q} -
q)\sqrt{1\over x},\cr a^{2m } &= 1 - \epsilon^{4m}{q^2\over x}, \,~~~~~
b^{2m} = 1 - \epsilon^{4m}{1\over x},
 ~~~~~~~\,  c^{2m} = \epsilon^{2m}\sqrt{1\over x}({1\over q} - q),~~~ ~~
d^{2m} = 0,\cr}
$$ and
$$ a+d ~: ~a-d~ :~ b+c ~: ~b-c  = {dn(u+\rho)\over dn(u-\rho) } : 1 :
{cn(u+\rho)\over cn(u-\rho)} : {sn(u+\rho)\over sn(u-\rho)}.
$$

\vskip1.5cm

\no {\bf Acknowledgements.} I thank Moshe Flato,   Kajiwara and Tetsuji
Miwa for discussions, and H. Araki  for generous hospitality at the
Research Institute for Mathematical Sciences, Kyoto University.  I thank
the Ministry of Education, Science, Sports and Culture for financial
support.

\ve

\vskip1in
\no {\bf References.}

\no [BD]  A.A. Belavin and V.G. Drinfeld, Sov.Sci.Rev.Math.{\bf 4} (1984)
93-165.
 
\no [F] ~ ~C. Fr\o nsdal, "Generalization and Deformations of Quantum
Groups;  Quantization of \line {\quad  
\quad ~all simple Lie Bi-algebras." q-alg 951000\hfill} 

\no [J] ~~~M. Jimbo, Commun.Math.Phys. 
 
\no [K] ~~V.G. Kac, {\it Infinite Dimensional Lie Algebras}, Cambridge
University Press 1990.

\no [L]~~~G. Lusztig, {\it Quantum Groups} 1993.

\no [R]~~ N.Yu. Reshetikhin, Lett.Math.Phys. {\bf 20} (1990) 331-336.

\end